\documentclass[12pt]{article}
\usepackage{amssymb,amsthm,amsmath,amsfonts,latexsym,tikz,hyperref,ytableau,authblk}
\usepackage[hmargin=1in,vmargin=1in]{geometry}

\newcommand{\ben}{\begin{enumerate}}
\newcommand{\een}{\end{enumerate}}
\newcommand{\ble}{\begin{lem}}
\newcommand{\ele}{\end{lem}}
\newcommand{\bth}{\begin{thm}}
\renewcommand{\eth}{\end{thm}}
\newcommand{\bpr}{\begin{prop}}
\newcommand{\epr}{\end{prop}}
\newcommand{\bco}{\begin{cor}}
\newcommand{\eco}{\end{cor}}
\newcommand{\bcon}{\begin{conj}}
\newcommand{\econ}{\end{conj}}
\newcommand{\bde}{\begin{defn}}
\newcommand{\ede}{\end{defn}}
\newcommand{\bex}{\begin{exa}}
\newcommand{\eex}{\end{exa}}
\newcommand{\barr}{\begin{array}}
\newcommand{\earr}{\end{array}}
\newcommand{\btab}{\begin{tabular}}
\newcommand{\etab}{\end{tabular}}
\newcommand{\beq}{\begin{equation}}
\newcommand{\eeq}{\end{equation}}
\newcommand{\bea}{\begin{eqnarray*}}
\newcommand{\eea}{\end{eqnarray*}}
\newcommand{\bal}{\begin{align*}}
\newcommand{\bce}{\begin{center}}
\newcommand{\ece}{\end{center}}
\newcommand{\bpi}{\begin{picture}}
\newcommand{\epi}{\end{picture}}
\newcommand{\bpp}{\begin{picture}}
\newcommand{\epp}{\end{picture}}
\newcommand{\bfi}{\begin{figure} \begin{center}}
\newcommand{\efi}{\end{center} \end{figure}}
\newcommand{\bprf}{\begin{proof}}
\newcommand{\eprf}{\end{proof}\medskip}

\newcommand{\bsl}{\begin{slide}{}}
\newcommand{\esl}{\end{slide}}
\newcommand{\bfr}{\begin{frame}}
\newcommand{\efr}{\end{frame}}

\newcommand{\hqed}{\hfill \qed}

\newcommand{\eqqed}[1]{$\rule{1ex}{0ex}\hfill{\dil#1}\hfill\qed$}

\newcommand{\ol}{\overline}

\newcommand{\hso}[1]{\hspace{-1pt}}
\newcommand{\vs}[1]{\vspace{#1}}

\newcommand{\emp}{\emptyset}

\newcommand{\sbe}{\subseteq}

\newcommand{\ptn}{\vdash}

\newcommand{\fl}[1]{\lfloor #1 \rfloor}

\def\<{\langle}
\def\>{\rangle}

\newcommand{\de}{\delta}

\newcommand{\la}{\lambda}

\newcommand{\bbN}{{\mathbb N}}
\newcommand{\bbP}{{\mathbb P}}

\newcommand{\cB}{{\cal B}}

\newcommand{\fS}{{\mathfrak S}}

\DeclareMathOperator{\des}{des}

\newcommand{\dil}{\displaystyle}

\newtheorem{thm}{Theorem}[section]
\newtheorem{prop}[thm]{Proposition}
\newtheorem{cor}[thm]{Corollary}
\newtheorem{lem}[thm]{Lemma}
\newtheorem{conj}[thm]{Conjecture}
\newtheorem{exa}[thm]{Example}

\DeclareMathOperator{\std}{std}

\newcommand{\fP}{{\mathfrak P}}
\newcommand{\ub}{\ol{u}}
\newcommand{\wb}{\ol{w}}

\DeclareMathOperator{\lwi}{lwi}

\begin{document}
\pagestyle{plain}

\title{Properties of  plactic monoid centralizers
}
\author[1]{Bruce E. Sagan}
\author[2]{Chenchen Zhao}
\affil[1]{Department of Mathematics, Michigan State University, East Lansing, MI 48824}
\affil[2]{Department of Mathematics, University of California, Davis, CA 95616}

\date{\today\\[10pt]
	\begin{flushleft}
	                                       \\[5pt]
	\end{flushleft}}

\maketitle

\begin{abstract}
Let $u$ be a word over the positive integers $\bbP$.  Motivated by a question involving crystal graphs, Sagan and Wilson initiated the study of the centralizer of $u$ in the plactic monoid which is the set
$$
C(u) = \{w \mid \text{$uw$ is Knuth equivalent to $wu$}\}.
$$
In particular, they conjectured the following stability phenomenon:   for any $u$ there is a positive integer $K$ depending only on $u$ such that $C(u^k) = C(u^K)$ for $k\ge K$.  We prove that this property holds for various $u$ including  words consisting of only ones and twos, as well as  permutations.
Sagan and Wilson also considered $c_{n,m}(u)$  which is the number of $w\in C(u)$ of length $n$ and maximum at most $m$.  They showed that 
$c_{n,m}(1)$ is a polynomial in $m$ of degree $n-1$ and conjectured  properties of the coefficients when it is expanded in a binomial coefficient basis.  We prove  some of these conjectures, for example, that the coefficients are always nonnegative integers.
\end{abstract}

\section{Introduction}

We will use the notation $\bbP$ and $\bbN$ for the positive and nonnegative integers, respectively.  And if $n\in\bbN$, then we let
$$
[n]=\{1,2,\ldots,n\}.
$$
Given a finite set $S$ we will use either $\#S$ or $|S|$ for its cardinality.  We will also use this notation for words $w$ over $S$, and $|w|$ will be called the {\em length} of $w$.
The {\em Kleene closure} of $S$, denoted $S^*$, is the set of all words which can be formed from elements of $S$.

We assume that the reader is familiar with the basic properties of the Robinson-Schensted-Knuth (RSK) bijection as well as the jeu-de-taquin (jdt) of Sch\"utzenberger.  The necessary background can be found in the texts of Sagan~\cite{sag:sym,sag:aoc} or Stanley~\cite{sta:ec2}.  We will let $P(w)$ denote the insertion tableau of $w$ under RSK.  And we will write $v\equiv w$ if $v$ and $w$ are Knuth equivalent, that is, $w$ can be obtained from $v$ by a sequence of Knuth transpositions.  As proved by Knuth~\cite{knu:pmg},
$$
v\equiv w \text{ if and only if } P(v)=P(w).
$$
Lacoux and Sch\"utzenberger's  {\em plactic monoid}~\cite{LS:mp} is $\bbP^*$ modulo Knuth equivalence equipped with the operation of concatenation.

Motivated by a question in representation theory, Sagan and Wilson~\cite{SW:cpm} studied the {\em centralizer} of a word $u$ in the plactic monoid which can be described as
\begin{align*}
    C(u) &= \{w \mid uw \equiv wu\}\\
    &=\{w \mid P(uw)=P(wu)\},
\end{align*}
where $uw$ denotes the concatenation of $u$ with $w$.
For example, they characterized $C(u)$ for various $u$ and studied an enumerative invariant of the centralizer, $c_{n,m}(u)$, defined below.  The purpose of the present work is to resolve or make progress on various conjectures that they posed in this context.

Our first results have to do with a stability phenomenon.  Call $u$ {\em stable} if there is a constant $K$ depending only on $u$ such that
$$
C(u^k) = C(u^K) \text{ for } k\ge K,
$$
where $u^k$ is $u$ concatenated with itself  $k$ times.
If we wish to be specific about the bound then we will write that $u$ is {\em $K$-stable}.
In particular, we say that $u$ is {\em strongly stable} if $K=1$.  We will establish some cases of the following conjecture.
\bcon[{\cite[Conjecture 7.2]{SW:cpm}}]
Every $u\in \bbP^*$ is  stable.
\econ
In the next section, we show that if $u$ consists of only ones and twos then it is strongly stable.  Then, in Section~\ref{per}, we show that any permutation is stable, and that the longest permutation is in fact strongly stable.

Sagan and Wilson also considered the cardinalities
$$
c_{n,m}(u) = \#\{w\in C(u) \mid |w|=n \text{ and } \max w\le m\}.
$$
In particular, they showed that $c_{n,m}(1)$ is a polynomial in $m$ of degree 
$n-1$.  In Theorem~\ref{cnm(1):th} below we give various properties of the coefficients of this polynomial when expanded in the basis of binomial coefficients $\binom{m}{k}$, $k\ge0$.  Parts (a), (b), and (e) of that result were conjectured by Sagan and Wilson.  We also show that the coefficients in the expansion of $c_{n,m}(1)$ in the basis $\binom{m-1}{k}$, $k\ge0$, have a nice combinatorial interpretation.

Various conjectures are scattered throughout.


\section{Words with only ones and twos}
\label{wot}

In this section we prove that if $u$ consists only of ones and twos then it is strongly stable.  In so doing, we characterize the tableaux $P(w)$ for $w\in C(u)$ when $u$ contains at least a single one and a single two. A characterization of $C(a^n)$ where $a\in\bbP$ and $n\ge1$ was given in~\cite{SW:cpm} and so, in particular, applies to any word which consists only of ones or only of twos.

We start with the case where $u=a^n$ for some $n\ge1$.  We will need the following previous result.

\bpr[{\cite[Theorem 5.2]{SW:cpm}}] 
\label{C(a^k)}
If  $a\in\bbP$, then $u=a$ is strongly stable.\hqed
\epr

\bco
\label{a^n}
If $a\in\bbP$ and $n\ge1$, then $u=a^n$ is strongly stable
\eco
\bprf
Suppose $k\ge1$. Applying the previous proposition twice
$$
C(u^k) = C(a^{nk}) = C(a) = C(a^n) = C(u)
$$
as desired.
\eprf

To deal with the case where $u$ contains both ones and twos, we need some preliminary lemmas about insertion tableaux.
If $P$ is a semistandard Young tableau then we let
$$
R_i(P) = \text{ the $i$th row of $P$}.
$$
We set $R_i(P)=\emp$ if $i$ is larger than the number of rows of $P$.
Similarly, if $u$ is a word then define
$$
R_i(u) = R_i(P(u)).
$$
A {\em singleton column} of $P$ is a column $C$ with $\#C=1$.  And we call $C$ a {\em singleton $a$-column} if it is a singleton column whose sole entry is $a$.
Finally, the {\em restriction of $w$ to ones and twos} is the word $\wb$ obtained by removing any elements of $w$ which are greater than or equal to three and keeping the same relative order of the remaining elements.  For example, the restriction of $w=3122413321$ is $\wb=122121$.  We define the {\em restriction} of an SSYT to ones and twos similarly.
\begin{lem}\label{lem:w'}
Let $u\in[2]^n$ and suppose $w$ has $P(w)$ with at least two rows and satisfying
$$
\max R_i(w) \le 2 \text{ for $i\le2$}.
$$ 
Then
\ben
\item[(a)] For $i\le 2$ we have $R_i(wu)=R_i(\wb u)$ and $R_i(uw)=R_i(u\wb)$ where $w'$ is the restriction of $w$ to its ones and twos.
\item[(b)] For $i\ge 3$ we have $R_i(wu)=R_i(w)=R_i(uw)$.
\een
\end{lem}
\bprf
First consider $wu$.  We compute $P(wu)$ by RSK.  So, first $P(w)$ is computed.  Since the first two rows of $P(w)$ consist only of ones and twos, and such elements are never bumped by numbers greater than two, we have $R_i(w)=R_i(\wb)$ for $i=1,2$.  Since $u$ consists only of ones and twos, and these can only occur in the first two rows of an insertion tableau, the previous sentence implies that $R_i(wu)=R_i(\wb u)$ for $i\le2$.  Also, since the first two rows of $P(w)$ do not contain any elements greater than two, the insertion of $u$ into $P(w)$ will not bump any elements from the second row into the third.  Thus $R_i(wu)=R_i(w)$ for $i\ge 3$, completing this case.

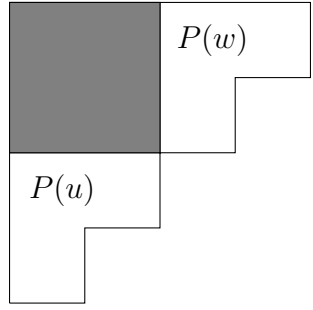
\begin{figure}
    \centering
\begin{tikzpicture}
\draw (0,0)--(0,4)--(4,4)--(4,3)--(3,3)--(3,2)--(2,2)--(2,1)--(1,1)--(1,0)--(0,0) (0,2)--(2,2)--(2,4); 
\draw(.7,1.5) node{$P(u)$};
\draw(2.7,3.5) node{$P(w)$};
\draw[fill=gray] (0,2) rectangle (2,4);
\end{tikzpicture}
    \caption{Computing $P(uw)$ via jdt}
    \label{jdt:fig}
\end{figure}

Now compute $P(uw)$ using jdt, starting with $P(u)$ southwest of $P(w)$  as in Figure~\ref{jdt:fig} and filling the skewing rectangle of cells row by row, right to left.  We claim that if $P_i$ is the skew tableau obtained after filling row $i$, then, for $i\ge 2$, we have the concatenations 
\begin{align*}
R_i(P_i) &=R_1(u)\ R_i(w),\\
R_{i+1}(P_i) &= R_2(u)\ R_{i+1}(w),\\
R_j(P_i) &= R_j(w) \text{ for $j\ge i+2$.}
\end{align*}
We prove this by reverse induction on $i$, where the base case before any slide as been applied is trivial.  Suppose these equations are true for $i+1$ and consider filling the cells of row $i$.  Since $i\ge 2$ we have $\max P(u)\le 2\le \min R_i(w)$.  So the first slide into row $i$ will bring the last element of row $1$ of $P(u)$ up.  If that element is in a singleton column of $P(u)$, then the slide will continue to move $R_{i+1}(w)$ one cell to the left.  Note that no elements of $R_{i+2}(w)$ will be affected since, by the induction assumption and the fact that in this case the second row of $P(u)$ is shorter than the first, row $i+2$ of $P(w)$ is already at least one cell to the left of row $i+1$.  Once all elements in singleton columns have been raised, applying jdt to a doubleton column will merely raise both elements.  This shows that the desired descriptions of the rows  still hold.

Now consider filling the first row.  We have shown that the equalities above hold when $i=2$.  Using arguments similar to the previous paragraph, one sees that once the first row is filled we have $R_3(uw)=R_3(w)$.  And we have previously shown $R_j(uw)=R_j(w)$ for $j\ge i+2=4$.  This verifies part (b) of the lemma for $uw$.  Furthermore, since the first two rows of $P(w)$ contain no elements which are three or larger, the restriction of $P_2$ to its ones and twos is exactly the same as the jdt tableau for computing $P(u\wb)$ after the row above $P(u)$ has been filled.
This imples that filling the first row of $P_2$ will cause exactly the same slides as filling the first row in computing $P(u\wb)$.  Part (a) of the lemma now follows.
\eprf

To prove our next lemma, we will need a couple of results from~\cite{SW:cpm}.
\bpr[{\cite[Lemma 3.3]{SW:cpm}}]
\label{row_cond}
Given $u$ and $w \in C(u)$ we let $P = P(w)$ have rows $R_i$ for $i\ge 1$. Also let $m = \max u$. If $u$ contains a subsequence $m,m-1,\dots,m-k+1$, then 
\[\max R_i \le m\] 
for $1\le i\le k$. \hqed
\epr

To state the second result, we need the notation
$$
\lwi(w) = \text{the length of a longest weakly increasing subsequence of $w$},
$$
and, for $a\in\bbP$,
$$
\lwi(w,a) = \text{the length of a longest weakly increasing subsequence of $w$ ending in $a$}.
$$

\bpr[{\cite[Corollary 4.3]{SW:cpm}}] 
\label{C(1)}
The following are equivalent.\ben
\item [(a)] $w \in C(1).$
\item [(b)] The entries of $R_1(w)$ are all ones.
\item [(c)] $\lwi(w) = \lwi(w,1)$. \hqed
\een
\epr

The following notation will be very useful.
$$
m_a(u) = \text{ the multiplicity of $a\in\bbP$ in the word $u$}.
$$

\begin{lem} \label{lem:(b)} Let $u \in [2]^n$ satisfy $m_1(u),m_2(u)\ge1$. If $w \in C(u)$, then $\max R_i(w) \le 2$ for $i \le 2$.
\end{lem}
\bprf
Suppose $w \in C(u)$ so that $P(uw) = P(wu)$. If $u$ contains the subsequence 21, then by Proposition~\ref{row_cond} with $m = 2$ and $k = 2$, we have $\max R_i(w) \le 2$ for $i\le 2$. Hence we may assume that $u = 1^k2^{n-k}$ with $1 \le k\le n-1$. By Proposition~\ref{row_cond} with $m = 2$ and $k = 1$, we have $\max R_1(w) \le 2$. To show that $\max R_2(w) \le 2$, we argue by contradiction. 

Suppose that $R_2(w)$ contains some entry greater than 2. If $R_1(w)$ contains a 2, then inserting a 1 from $u$ will bump the leftmost 2 from the first row, which in turn will bump the leftmost entry greater than 2 from the second row. If follows from~\cite[Corollary 3.2]{SW:cpm} that inserting $u$ into $P(w)$ cannot bump any element not in $u$, so this is a contradiction. Hence, the entries of $R_1(w)$ must be all 1's. By Propositions~\ref{C(a^k)} and~\ref{C(1)}, this implies $w \in C(1) = C(1^k)$.  And $u\in C(w)$.  Together these observations yield
$$
 P(1^k2^{n-k}w)=P(w1^k2^{n-k}) = P(1^kw2^{n-k}).
$$

Now consider $P(1^k v)$ for any word $v$.  Using jdt as in Lemma~\ref{lem:w'}, we see that we have the concatenation
$R_1(1^k v) = 1^k R_1(v)$ as well as $R_i(1^k v) =R_i(v)$ for $i\ge2$. Combining this with the previously displayed equations shows that $P(2^{n-k} w) = P(w 2^{n-k})$.
So, $w \in C(2^{n-k}) = C(2)$ by Proposition~\ref{C(a^k)}. By~\cite[Theorem 4.2]{SW:cpm}, every column of $P(w)$ must contain a 2, forcing $R_2(w)$ to consist entirely of 2's since $R_1$ consists entirely of 1's. This contradicts our assumption that $R_2(w)$ contains an entry greater than 2. Thus, $\max R_i(w) \le 2$ for $i\le2.$
\eprf

To state a more precise relationship between the second rows of $P(wu)$, $P(w)$, and $P(u)$ when $u,w$ are all ones and twos we will need the notation
$$
c_a(u)=\text{ number of singleton $a$-columns in $P(u)$},
$$
for $a\in\bbP$.

\ble
\label{R2(wu)}
If $u,w\in[2]^*$, then
$$
\#R_2(wu) = \#R_2(w) + \#R_2(u) +\min(c_1(u), c_2(w)).
$$
\ele
\bprf
To compute $P(wu)$, we can take any word Knuth equivalent to $u$ and insert it in $P(w)$.  So we can assume that $u$ decomposes as a concatenation of 
$$
u' = (21)^{\#R_2(u)} \text{ and } u'' = 1^{c_1(u)}2^{c_2(u)}.
$$
By construction, $P(u) = P(u'u'')$ and hence $P(wu)=P(wu'u'')$. The net result of inserting the word $u'$ into $P(w)$ appends each 1 from $u'$ to the first row and each 2 from $u'$ to the second row. Thus, we have $R_i(wu') = i^{\#R_2(u)} \ R_i(w)$ and  $c_i(wu') = c_i(w)$ for $i\le 2$. Next, we insert $u''$ into $P(wu')$. During this insertion, every 1 from $u''$ bumps a 2 from the first row whenever a 2 is present, and each bump contributes a 2 to the second row. Thus, the number of 2's bumped off into the second row of $P(wu')$ is 
$$
\min(m_1(u''),c_2(wu')) = \min(c_1(u),c_2(w)).
$$
Hence, the total number of 2's in $R_2(wu)$ is 
$$
\#R_2(wu) = \#R_2(w) + \#R_2(u) +\min(c_1(u), c_2(w)).
$$
as desired.
\eprf

We will now give characterizations of words in the centralizer of a $u$ consists of only ones and twos.  There will be  three cases depending on the relative multiplicities of $1$ and $2$ in $u$.

\bth\label{thm <} Let $u \in [2]^n$ with $1\le m_1(u) < m_2(u)$. Then $w \in C(u)$ if and only if $P(w)$ satisfy following two conditions.
\begin{enumerate}
\item [(a)] $c_1(w) = c_1(u) \le c_2(w)$ or $c_1(w) = c_2(w)< c_1(u)$.
\item[(b)] $\max R_i(w) \le 2$ for $i \le 2$.
\end{enumerate}
\eth

\bprf
We first show that (a) and (b) imply $w \in C(u).$ Let $\wb$ be the restriction of $w$ to its 1s and 2s. Since (b) holds, we can apply Lemma~\ref{lem:w'} to obtain $R_i(wu) = R_i(\wb u)$ and $R_i(uw) = R_i(u\wb)$ for $i \le 2$ and $R_i(wu) = R_i(uw)$ for $i \ge 3.$ Hence, it suffices to verify that $R_i(\wb u) = R_i(u\wb)$ for $i \le 2$.  But since $\wb$  and $u$ only consist of 1s and 2s this is equivalent to proving $P(\wb u) = P(u\wb).$

Applying Lemma~\ref{R2(wu)} we see that the total number of 2's in $R_2(\wb u)$ is 
\beq
\label{R2(w'u)}
\#R_2(\wb u) = \#R_2(\wb) + \#R_2(u) +\min(c_1(u), c_2(\wb)).
\eeq
Interchanging the roles of of $\wb$ and $u$ in the same lemma gives
\beq
\label{R2(uw')}
\#R_2(u\wb) =\#R_2(u) + \#R_2(\wb) +\min(c_1(\wb), c_2(u)).
\eeq
Since $c_i(\wb) = c_i(w)$ for $i\le 2$ and  the given condition $m_1(u) <m_2(u)$ implies $c_1(u)<c_2(u)$, we now compare the two minima using the restrictions in part (a) of the theorem. If $c_1(u) = c_1(w) \le c_2(w)$, then 

$$
\min(c_1(u),c_2(\wb)) =c_1(u) = c_1(w) = \min(c_1(\wb),c_2(u)).
$$
If $c_1(u) > c_1(w) = c_2(w)$, then
$$
\min(c_1(u),c_2(\wb)) = c_2(w) = c_1(w) = \min(c_1(\wb),c_2(u)).
$$
Thus, $\#R_2(\wb u) = \#R_2(u\wb)$. Since $\wb u$ and $u\wb$ contain the same number of 1s and 2s and no other elements, their $P$-tableaux are uniquely determined by the length of the second row.  It follows that  $P(\wb u) = P(u\wb)$ so that, by the previous discussion,  $w \in C(u).$

For the converse, suppose that $w \in C(u)$, so that $P(wu) = P(uw)$. Since $u$ contains both one and two,  condition (b) holds by Lemma~\ref{lem:(b)}. Applying Lemma~\ref{R2(wu)} again to both $wu$ and $uw$, the equality $R_2(wu) = R_2(uw)$ yields 
$$
\#R_2(w) + \#R_2(u) +\min(c_1(u), c_2(w)) = \#R_2(u) + \#R_2(w) +\min(c_1(w), c_2(u)),
$$
which simplifies to 
\beq
\label{min(c1,c2)}
\min(c_1(u),c_2(w)) = \min(c_1(w),c_2(u)).
\eeq
Since $m_2(u) > m_1(u)$, we have $c_2(u) > c_1(u)$. In particular, 
$$
c_2(u) > \min(c_1(u),c_2(w)) = \min(c_1(w),c_2(u)).
$$
Thus, both minima must equal $c_1(w)$, which implies that either $c_1(w) = c_1(u) \le c_2(w)$ or $c_1(w) = c_2(w) \le c_1(u).$
\eprf

\bth\label{thm =} Let $u \in [2]^n$ with $ m_1(u) = m_2(u)\ge1$. Then $w \in C(u)$ if and only if $P(w)$ satisfy the following two conditions.
\begin{enumerate}
    \item[(a)] $\min(c_1(w), c_2(w)) \ge c_1(u)$ or $c_1(w) = c_2(w) < c_1(u).$
    \item[(b)] $\max R_i(w) \le 2$ for $i \le 2$.
\end{enumerate}
\eth

\bprf The proof follows the same structure as the demonstration of the previous theorem. We first show that (a) and (b) imply $w \in C(u)$. As before, let $\wb$ be the restriction of $w$ to its ones and twos so that it suffices to show  $P(\wb u) = P(u\wb)$.

We compare $\#R_2(\wb u)$ and $\#R_2(u\wb)$ using equations~\eqref{R2(w'u)} and~\eqref{R2(uw')}. Now $c_i(\wb) = c_i(w)$ for $i\le 2$ and the given condition $m_1(u) =m_2(u)$ implies $c_1(u)=c_2(u)$. Now compare the two minima using the restrictions in part (a) of the theorem. If
$\min\{c_1(w), c_2(w)\} \ge c_1(u)$, then $\min\{c_1(\wb), c_2(\wb)\} \ge c_1(u)$ and
$$
\min(c_1(u),c_2(\wb)) =c_1(u) =  \min(c_1(\wb),c_1(u)) = \min(c_1(\wb),c_2(u)).
$$
If $c_1(w) = c_2(w) < c_1(u)$, then
$$
\min(c_1(u),c_2(\wb)) = c_2(w) = c_1(w) = \min(c_1(\wb),c_2(u)).
$$
So, in either case, $\#R_2(\wb u) = \#R_2(u\wb)$. 
It follows from the same reasoning as in the proof of the previous theorem that $P(\wb u) = P(u\wb)$, which implies $w \in C(u).$

For the converse, suppose that $w \in C(u)$, so that $P(wu) = P(uw)$. 
As in the previous proof, condition (b) holds and we can assume equation~\eqref{min(c1,c2)}.
 Since $m_1(u) = m_2(u)$, we have $c_1(u) = c_2(u)$.
So~\eqref{min(c1,c2)} becomes
$$
 \min(c_1(u),c_2(w))=\min(c_1(u),c_1(w)).
$$
 
There are two cases depending on whether $c_2(w)\ge c_1(u)$ or $c_2(w)< c_1(u)$. In the first case, both minima are equal to $c_1(u)$, which implies $c_1(w)\ge c_1(u)$. Thus $\min(c_1(w), c_2(w)) \ge c_1(u)$. In the second case, the equality of minima gives $c_2(w) = \min(c_1(u),c_1(w))$. Since $c_1(u) > c_2(w)$, we must have $c_1(w) = c_2(w)  <c_1(u).$
\eprf

\bth \label{thm >}Let $u \in [2]^n$ with $1\le m_2(u)<m_1(u)$. Then $w \in C(u)$ if and only if $P(w)$ satisfy the following two conditions.
\begin{enumerate}
  \item[(a)] $c_2(w) = c_2(u) \le c_1(w)$ or $c_1(w) = c_2(w) < c_2(u)$.
  \item[(b)]  $\max R_i(w) \le 2$ for $i \le 2$.
\end{enumerate}
\eth

\bprf 
The proof again follows the same structure as the demonstration of the previous theorem. We first verify the reverse implication. Let $\wb$ be the restriction of $w$ to its ones and twos.  As before, it suffices to show $P(\wb u)=P(u\wb)$ which will be implied by proving  $\#R_2(\wb u) = \#R_2(u\wb)$.  And, by equations~\eqref{R2(w'u)} and~\eqref{R2(uw')}, this reduces to showing
$\min(c_1(u),c_2(\wb)) = \min(c_1(\wb),c_2(u))$.
We have $c_i(\wb)=c_i(w)$ for $i\le2$ and $m_1(u)>m_2(u)$ implies $c_1(u)>c_2(u)$.
{If the first condition in (a) holds then $c_2(\wb) = c_2(u) \le c_1(\wb)$}.  So
$$
\min(c_1(u),c_2(\wb)) =c_2(\wb) = c_2(u) = \min(c_1(\wb),c_2(u)).
$$
On the other hand, the second condition implies $c_1(\wb) = c_2(\wb) < c_2(u)$.  Thus
$$
\min(c_1(u),c_2(\wb)) = c_2(\wb) = c_1(\wb) = \min(c_1(\wb),c_2(u)).
$$
So in either case we have equal minima.

For the converse, suppose that $w \in C(u)$, so that $P(wu) = P(uw)$. Then condition (b) holds by Lemma~\ref{lem:(b)}. And, as before, equation~\eqref{min(c1,c2)} holds.
Since $m_1(u) > m_2(u)$, we have $c_1(u) > c_2(u)$, and hence \[c_1(u)> \min(c_1(w),c_2(u)) = \min(c_1(u),c_2(w)). \]Both minima must therefore equal $c_2(w)$, which implies that either $c_2(w) = c_2(u) \le c_1(w)$ or $c_1(w) = c_2(w) < c_2(u)$.
\eprf

Using these characterizations of words $u$ over $[2]$, we can now prove the strong stability of $u$.
First we need two final lemmas comparing the number of singleton columns in $P(u)$ and $P(u^k)$.

\ble\label{lem:c1} Let $u \in [2]^n$ with $1 \le m_1(u) \le m_2(u)$. For all $k\ge1$ we have  $c_1(u^k) = c_1(u)$.
\ele

\bprf The proof is by induction on $k$, where the  case $k=1$ is obvious. Assume $c_1(u^k) = c_1(u)$ for some $k \ge 1$ and consider $u^{k+1} = u^k u.$ Applying Lemma \ref{R2(wu)}, inserting $u$ into $P(u^k)$ gives 
\begin{align*}
    \#R_2(u^{k+1}) &= \#R_2(u^k) + \#R_2(u) + \min(c_1(u^k), c_2(u)).
\end{align*}
Using the inductive hypothesis we obtain $c_1(u^k)=c_1(u)$.  And the assumption $m_1(u)\le m_2(u)$ implies $c_1(u)\le c_2(u)$.  Thus
\[
\#R_2(u^{k+1}) = \#R_2(u^k) + \#R_2(u) + c_1(u).
\]

For any word $v$ consisting only of ones and twos, its insertion tableau $P(v)$ has two rows, and
\[
c_1(v) = m_1(v) - \#R_2(v).
\]
Since $m_1(u^{k+1}) = (k+1)m_1(u)$, substitution into the previous equation gives
\begin{align*}
c_1(u^{k+1})
&= (k+1)m_1(u) - \bigl(\#R_2(u^k) + \#R_2(u) + c_1(u)\bigr)\\
&= \bigl(km_1(u) - \#R_2(u^k)\bigr) + \bigl(m_1(u) - \#R_2(u) - c_1(u)\bigr)\\
&= c_1(u^k),
\end{align*}
where the last equality follows from the identities 
$c_1(u^k)=k m_1(u)-\#R_2(u^k)$ and $c_1(u)=m_1(u)-\#R_2(u)$. Thus $c_1(u^{k+1})=c_1(u^k) = c_1(u)$ which completes the induction.
\eprf 

\ble\label{lem:c2} Let $u \in [2]^n$ with $1 \le m_2(u)  \le m_1(u)$. For all $k\ge1$ we have  $c_2(u^k) = c_2(u)$.
\ele

\bprf
The argument follows the same inductive pattern as in Proposition~\ref{lem:c1}.  {In particular, in the induction step one uses Lemma~\ref{R2(wu)} together with} $m_2(u)\le m_1(u)$, so that $\min(c_1(u),c_2(u))=c_2(u)$ yields
\[
\#R_2(u^{k+1}) = \#R_2(u^k) + \#R_2(u) + c_2(u).
\]
Since $c_2(v)=m_2(v)-\#R_2(v)$ for binary $v$, {an analogous} computation as before gives 
$c_2(u^{k+1})=c_2(u^k)=c_2(u)$, completing the induction.
\eprf

We now use our characterizations to prove strong stability.

\bth 
\label{[2]^nStrong}
If $u \in [2]^n$  then $u$ is strongly stable.
\eth
\bprf 
If $m_1(u)=0$ or $m_2(u)=0$ then this is a special case of Corollary~\ref{a^n}.

The other possibility is that $m_1(u),m_2(u)\ge1$.
Clearly $m_i(u^k) = km_i(u)$ for $i = 1,2$. Hence, $u$ and $u^k$ both satisfy one of $m_1>m_2$, $m_1=m_2$, or $m_1<m_2$. Thus, both $u$ and $u^k$ fall under the same hypothesis of Theorems~\ref{thm <},~\ref{thm =}, {or}~\ref{thm >}, respectively.

In each of these three theorems, condition (b) is identical for $u$ and $u^k$ since it does not involve $u$ at all.
For condition~(a), Lemmas~\ref{lem:c1} and~\ref{lem:c2} show that the relevant parameters remain unchanged upon repetition, that is
\begin{align*}
c_1(u^k) &= c_1(u) \text{ if  $m_1(u)\le m_2(u)$, and}\\
c_2(u^k) &= c_2(u) \text{ if $m_1(u)>m_2(u)$}.
\end{align*}

Thus, condition~(a) is also identical for $u$ and $u^k$ {in all cases}.
Since both conditions (a) and (b) agree, the theorems show that  
{$w \in C(u)$ if and only if  $w \in C(u^k)$,}
and hence $C(u) = C(u^k)$.\eprf

\section{Permutations}
\label{per}

Let $\fS_m$ be the symmetric group of permutations $u$ of  $[m]$ with elements expressed in one-line form, that is, as rearrangements of $[m]$.  We wish to show that every such $u$ is $m$-stable.  And the decreasing permutation $\de_m=m(m-1)\ldots 1$  is, in fact, strongly stable.  We use this notation rather than $w_0$ (which is usually employed for the longest element of a Coxeter group) to show the dependence on $m$.
We first need some preliminary results.  

We will use a theorem of Greene~\cite{gre:est} about the shape of the $P$-tableau under RSK.  Call a sequence $x$ {\em weakly $i$-increasing} if $x$ can be written as a union of $i$ disjoint weakly increasing subsequences.  For $u\in\bbP^*$ we define
$$
\lwi_i(u)=\text{length of a longest weakly $i$-increasing subsequence of $u$.}
$$
\bth[\cite{gre:est}]
\label{greene}
If $u\in \bbP^*$ then let $P(u)$ have shape 
$\la=(\la_1,\la_2,\ldots,\la_l)$.  For any $i\ge 1$ we have

\vs{5pt}

\eqqed{
\la_1+\la_2+\cdots+\la_i = \lwi_i(u).
}
\eth

In order to use Greene's theorem, we will need the following result.

\ble
\label{lwi}
If $u\in\fS_m$ and $i,k\ge 1$ then
$$
\lwi_i(u^{k+1})\ge\lwi_i(u^k) + i.
$$
\ele
\bprf
Let $x$ be a longest $i$-increasing subsequence of $u^k$.  Write $x$ as a union of subsequences
$$
x=x^{(1)}\uplus \cdots \uplus x^{(i)}
$$
where each subsequence is increasing and $x^{(1)},\ldots,x^{(s)}$ are those which end in $m$ for some $0\le s\le i$.  Further, if $s\ge1$ then suppose that $x^{(1)}$ contains the leftmost of these $m$'s in $u^k$.  Since $u$ only contains one $m$, it must be that these $s$ subsequences contain $m$'s from different factors of $u$ in $u^k$.  It follows that the sequences
$$
y^{(1)}=x^{(1)}m^{s-1},
$$
and
$$
y^{(j)}=\text{$x^{(j)}$ with $m$ removed}
$$
for $2\le j\le s$ are disjoint increasing subsequences of $u^k$.  Furthermore
$$
y:= y^{(1)}\uplus \cdots y^{(s)}\uplus x^{(s+1)}\uplus\cdots\uplus x^{(i)}
$$
has $|y|=|x|$ and $y^{(1)}$ is the only subsequence which can end in an $m$.

Continue this process with the subsequence obtained from $y$ by disregarding $y^{(1)}$ and considering all the remaining subsequences ending in $m-1$, etc.
One finally obtains a subsequence 
$$
z=z^{(1)}\uplus \cdots \uplus z^{(i)}
$$
of $u^k$ with $|z|=|x|$ and $\max z^{(j)}\le m-j+1$ for $j\in[i]$.  Now we have a $i$-increasing subsequence in $u^{k+1}$ of the form
$$
z=z^{(1)}m \uplus z^{(2)}(m-1) \uplus\cdots \uplus z^{(i)}(m-i+1).
$$
Thus 
$$
\lwi_i(u^{k+1})\ge|z|=|x|+i = \lwi_i(u^k) + i
$$
as desired.
\eprf

We now show that, for $u\in\fS_m$ and large enough $k$, the rows of $P(u^{k+1})$ are similar to those of $P(u^k)$.  To see how that rows of these tableaux are related, we first need a lemma.  Recall that
$$
\de_m = m(m-1)\ldots 1.
$$
\ble
\label{de}
If $v\in[m]^*$ then for all $i\in[m]$
\beq
\label{Ridv}
R_i(\de_m v) = i R_i(v).
\eeq
\ele
\bprf
We have that $P(\de_m)$ is a single column with entries $1,2,\ldots,m$.  Now using jdt to compute $P(\de_m v)$ from $P(\de_m)$
and $P(v)$, we see that the elements of the column just move up one row each time an empty cell is filled.  This proves the result.
\eprf

It turns out that if $u\in\fS_m$ then $P(u^{k+1})$ and $P(u^k)$ are similarly related for $k\ge m$.
Generalizing from the previous section we define, for any $u\in\bbP^*$, the {\em restriction of $u$ to $[m]$} to be the subsequence $\ub$ of $u$ obtained by removing all elements of $u$ larger than $m$.

\ble Suppose that $u \in \fS_m$ and $k\ge m$. Then for all 
$i \in [m]$, we have 
\beq
\label{Riuk}
R_i(u^{k+1}) = iR_i(u^k).
\eeq
\ele
\bprf
We will induct on $m$ where the result is clear from $m=1$.
So, suppose $m\ge 2$ and 
let $\overline{u}$ denote the restriction of $u$ to $[m-1]$. 
By induction, we have 
\beq
\label{R1ubk}
R_i(\ub^{k+1}) = i R_i(\ub^k)
\eeq
for $i\in[m-1]$ and $k\ge m-1$.

We now induct on $k$.  In fact we claim that, by the previous lemma, it suffices to prove the base case when $k=m$.  To see this, assume that the current lemma holds for a given value of $k$.  Then, comparing equations~\eqref{Ridv} and~\eqref{Riuk}, we have the Knuth equivalence 
$u^{k+1}\equiv \de_m u^k$.  So,
$$
u^{k+2}=u^{k+1} u \equiv \de_m u^k u = \de_m u^{k+1}.
$$
Applying the previous lemma again gives $R_i(u^{k+2})=iR_i(u^{k+1})$ for $i\in[m]$ as desired.

To prove the lemma under consideration for $u$ when $k=m$, we will also induct on $i$.
First, consider the base case when $i=1$.  Note that by~\eqref{R1ubk} we
have
\beq
\label{R1ub}
R_1(\ub^m) = 1 R_1(\ub^{m-1}) \text{ and } 
R_1(\ub^{m+1}) = 1 R_1(\ub^m).
\eeq
Since $m$ cannot bump any smaller element we have the concatenation
\beq
\label{R1um}
R_1(u^m) = R_1(\overline{u}^m)m^a
\eeq
for some multiplicity $a\ge0$.  By similar reasoning and using~\eqref{R1ub}, we obtain
\beq
\label{R1um+1}
R_1(u^{m+1}) = R_1(\overline{u}^{m+1}) m^b = 1R_1(\overline{u}^m)m^b
\eeq
for some $b\ge0$.   To finish the base case, it suffices to show that $a = b$. 

Let $x$ denote a longest weakly increasing subsequence of $u^m$. Then, by Theorem~\ref{greene} and~\eqref{R1um}.
\beq
\label{R1m}|x| = |R_1(u^m)| = |R_1(\overline{u}^m)|+a.
\eeq 
Using~\eqref{R1um+1} and the fact that $xm$ is also a weakly increasing subsequence of $u^{m+1}$ we obtain
\beq
\label{R1m+1}
 1 + |R_1(\overline{u}^m)| + b=|R_1(u^{m+1})| \ge |xm| = |x|+1.
\eeq 
Comparing equations~\eqref{R1m} and~\eqref{R1m+1} and
we have
$$
b\ge a.
$$
On the other hand, the tableau $P(u^{m+1})$ can be obtained by inserting $u$ into $P(u^m)$, and $u$ has a single $m$.  It follows that
$$
b \le a+1.
$$
Comparing these two inequalities for $a,b$ we see that it suffices to show that some $m$ is displaced in passing from $R_1(u^m)$ to $R_1(u^{m+1}).$ We will divide the proof  into two cases based on the positivity of $a$.  In both cases we will use the notation
$$
j=|R_1(\ub^m)|.
$$

First suppose that  $a > 0$.  Then, by equation~\eqref{R1um} and the definition of $j$, we have an $m$ in cell $(1,j+1)$ of $R_1(u^m)$.
But from~\eqref{R1um+1} we see that cell $(1,j+1)$ of $R_1(u^{m+1})$ is occupied by an element of $\ub$.  So, the $m$ in this cell was bumped.

Now suppose that $a=0$ so that by~\eqref{R1um}
$$
|R_1(u^m)|=|R_1(\ub^m)|=j.
$$
But, from the first half of~\eqref{R1ub},
$$
|R_1(\ub^{m-1})|=j-1.
$$
It follows that the $m$ in the last copy of $u$ in $u^m$ must enter the $P$-tableau in cell $(1,j)$ and then be bumped from that cell by a later element.
Now comparison of the first and second equations in~\eqref{R1ub} show that the insertion of the $(m+1)$st copy of $\ub$ causes the same elements to be bumped as in the insertion of the $m$th copy, and from the same place but one cell to the right.  Thus the $m$ in the last copy of $u^{m+1}$ enters in cell $(1,j+1)$ and is bumped by a later element.  This completes the proof of the base case $i=1$.  The induction step is similar, so we will only give details for the differences.

Now suppose that equation~\eqref{Riuk} holds for values less than $i$.
Write
\beq
\label{Rium}
R_i(u^m) = R_i(\overline{u}^m)m^c
\eeq
and, using the induction on $m$,
\beq
\label{Rium+1}
R_i(u^{m+1}) = 1R_i(\overline{u}^m)m^d
\eeq
To prove $c=d$, first combine Theorem~\ref{greene} and Lemma~\ref{lwi} to obtain
$$
|R_1(u^{m+1})|+\cdots+|R_i(u^{m+1})|\ge |R_1(u^m)|+\cdots+|R_i(u^m)|+i.
$$
By induction on $i$ we have $|R_j(u^{m+1})|=|R_j(u^m)|+1$ for $j<i$.  Combining this with~\eqref{Rium} and~\eqref{Rium+1} gives $d\ge c$.
The induction on $i$ also shows that there is a single $m$ which is bumped out of row $i-1$ and thus into row $i$, giving $d\le c+1$.  There are now two cases depending on the positvity of $c$, but they are much like the base case and so the demonstration is omitted.
\eprf

Combining equations~\eqref{Ridv} and~\eqref{Riuk} we obtain the following.

\ble\label{Pduk} Suppose that $u \in \fS_m$. Then for all $k \ge m$, we have 

\hspace{10pt}

\eqqed{
P(u^{k+1}) = P(\de_m u^k).
}
\ele

We need one more result before we can prove the main theorem of this section

\bpr[{\cite[Theorem 5.3]{SW:cpm}}]
\label{staircase}
We have $w \in C(m(m-1)\dots1)$ if and only if $P= P(w)$ satisfies \[\max R_i \le m \text{ for all }1\le i \le m\] where $R_i$ is the $i$th row of $P$. \hqed
\epr

\bth If  $u \in \fS_m$ then $u$ is $m$-stable.
\eth

\bprf
It suffices to prove that if $k\ge m$ then
$$
C(u^k) = C(u^{k+1}).
$$
We first show that $C(u^k) \subseteq C(u^{k+1})$.  Suppose that $w \in C(u^k)$, so that 
\beq
\label{Pwuk}
wu^k \equiv u^k w.
\eeq
Let $P = P(w)$, and denote the $i$-th row of $P$ by $R_i$  for $i \ge 1$.  
Since {$k\ge m$ we have that} $u^k$ contains the decreasing subsequence $\de_m$ {obtained by picking $m-i+1$ from the $i$th factor of $u^k$.}
Proposition~\ref{row_cond} implies that 
\[
\max R_i \le m \qquad \text{for } 1 \le i \le m.
\]
Then, by Proposition~\ref{staircase}, we have $w \in C(\de_m)$.  
Applying Lemma~\ref{Pduk} twice along with the previous sentence and equation~\eqref{Pwuk}, we obtain the following sequence of 
{Knuth equivalences}
\[
wu^{k+1} \equiv w\de_m u^{k} \equiv \de_m wu^{k} \equiv \de_m u^kw \equiv u^{k+1}w,
\]
which shows that $w \in C(u^{k+1})$. 

Conversely, suppose that $w \in C(u^{k+1})$, so that 
\beq
\label{Pwuk+1}
P(wu^{k+1}) = P(u^{k+1}w).
\eeq
Again, by Propositions~\ref{row_cond} and~\ref{staircase}, we have 
$w \in C(\de_m)$, so that $P(\de_m w) = P(w\de_m)$.  
Combining this with  Lemma~\ref{Pduk} and equation~\eqref{Pwuk+1}, we can write
\[
{
P(\de_m wu^k) = P(w\de_m u^k) =  P(wu^{k+1})  = P(u^{k+1}w) = P(\de_m u^kw).
}
\]
The equality $P(\de_m wu^k) = P(\de_m u^kw)$ together with Lemma~\ref{de} imply
$P(wu^k) = P(u^k w)$.
Thus $w \in C(u^k)$ and $C(u^{k+1})\subseteq C(u^k)$.  This completes the proof of the theorem.
\eprf

We conjecture that the previous theorem can be generalized as follows.
Say that a word $u$ is {\em $m$-packed} if it contains at least one copy of every integer in $[m]$ where $m=\max u$.
This terminology was introduced by Pechenik~\cite{pec:csi} in regards to poset labelings used for a $K$-theoretic analogue of Sch\"utzenberger's promotion operator.  The orbit structure of $K$-promotion on $m$-packed labelings of rooted trees was shown to have interesting properties by Kimble, Sagan, and St.\ Dizier~\cite{KSS:kml}.  Obviously, every $u\in\fS_m$ is $m$-packed.
Computer calculations support the following conjecture about stability of $m$-packed words $u$.  Our computations were for all $u$ with $m=3,4$ of length at most $8$.  (The cases of $m=1$ and $2$  having already been proven in Corollary~\ref{a^n} and Theorem~\ref{[2]^nStrong}, respectively.)   Of course, we could not compute the whole set $C(u^k)$ since it is infinite.  But we calculated the set 
$C'(u^k)= C(u^k)\cap[m]^l$ for $l\le 10$ and 
verified that $C'(u^k)=C'(u^{k+1})$ for $m\le k\le 14$.
\begin{conj}
If $u$ is $m$-packed then $u$ is $m$-stable.
\end{conj} 

Note that elements of $\fS_m$ need not be strongly stable.
For example, take $u = 1234$. Then the word $w = 4123$ lies in $C(u^3)$ but not in $C(u^2)$.  However, the longest word is strongly stable.

\bpr
\label{C(u^k)} 
The permutation   $\de_m=m(m-1)\cdots 1$ is strongly stable.
\epr

\bprf
It was noted in~\cite[Lemma 5.1]{SW:cpm} that, for any $u\in\bbP^*$ and $k\ge1$, we have 
$C(u)\sbe C(u^k)$.  For the other containment,
suppose $w \in C(\de_m^k)$ and let $P = P(w)$ have rows $R_i$ for $i\ge 1.$ Since $\de_m^k$ contains the subsequence $m(m-1)\dots 1$, Proposition~\ref{row_cond} implies that $\max R_i \le m$ for $1\le i \le m.$ Thus $P$ satisfies the condition in Proposition~\ref{staircase}, so $w \in C(\de_m)$. 
\eprf


\section{Coefficients}

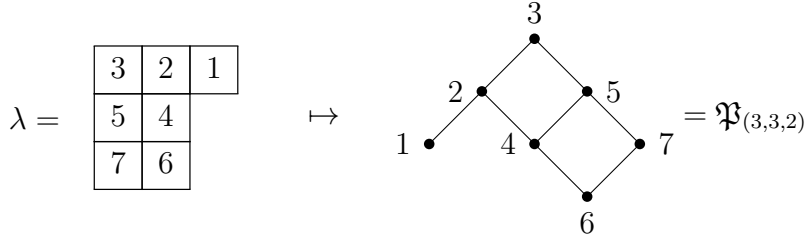
\begin{figure}
    \centering
\begin{tikzpicture}[scale=.7]
\draw(-7.5,1.5) node{$\la=$};
\draw(-5,1.5) node{
\begin{ytableau}
3&2&1\\
5&4\\
7&6
\end{ytableau}
};
\draw(-2,1.5) node{$\mapsto$};
\fill(3,0) circle(.1);
\fill(0,1) circle(.1);
\fill(2,1) circle(.1);
\fill(4,1) circle(.1);
\fill(1,2) circle(.1);
\fill(3,2) circle(.1);
\fill(2,3) circle(.1);
\draw (0,1)--(2,3)--(4,1)--(3,0)--(1,2) (2,1)--(3,2);
\draw(3,-.5) node{$6$};
\draw(-.5,1) node{$1$};
\draw(1.5,1) node{$4$};
\draw(4.5,1) node{$7$};
\draw(.5,2) node{$2$};
\draw(3.5,2) node{$5$};
\draw(2,3.5) node{$3$};
\draw(6,1.5) node{$=\fP_{(3,3,2)}$};
\end{tikzpicture}
    \caption{The poset $\fP_{(3,2,2)}$}
    \label{fP:fig}
\end{figure}

Let
$$
c_{n,m}(u) = \#\{w\in C(u) \mid |w|=n \text{ and } \max w\le m\}.
$$
The following result was proved in~\cite{SW:cpm} under the assumption the $n$ is fixed and $m\ge n$.  But an examination of the proof shows that the assumed inequality is not necessary for $n\ge2$.
\bpr[{\cite[Corollary 6.6(a)]{SW:cpm}}]
\label{cnm(1):pro}
Fix $n\ge2$ and let $m\ge0$ vary.  Then $c_{n,m}(1)$ is a polynomial in $m$ of degree $n-1$ with leading coefficient $1/(n-1)!$.
\epr
We can now prove part of a conjecture from~\cite{SW:cpm} about the the expansion of $c_{n,m}(1)$ in the binomial coefficient basis.  In the proof we will use the notation
$\la\ptn n$ if the integer partition $\la=(\la_1,\ldots,\la_l)$ satisfies $\sum_i \la_i=n$. The {\em Young diagram} of $\la$ consists of $l$ left-justified rows of boxes with $\la_i$ boxes in row $i$ from the top.  To each box we assign {\em coordinates} $(i,j)$ where $i$ and $j$ are the row and column numbers of the box.
A filling of the Young diagram of $\la=(3,2,2)$ is shown on the left in Figure~\ref{fP:fig} and the box containing $6$ has coordinates $(3,2)$.
We will associate with $\la$ a labeled  poset (partially ordered set) $\fP_\la$, defined as follows.  Label the boxes of the Young diagram of $\la$ with $[n]$ where $\la\ptn n$ by putting $1,\ldots,\la_1$ in the first row from right to left, then $\la_1+1,\ldots,\la_1+\la_2$ in the second row from right to left, and so forth.  
See the labeling of $(3,2,2)$ on the left in Figure~\ref{fP:fig} as an example.
We now construct the poset $\fP_\la$ on $[n]$ by letting
$$
\text{$k\le k'$ in  $\fP_\la$ if and only if $i'\ge i$ and $j'\ge j$
}
$$
where $k$ and $k'$ have coordinates $(i,j)$ and $(i',j')$, respectively, in the labeling of $\la$.  The poset $\fP_{(3,2,2)}$ is shown on the right in Figure~\ref{fP:fig}.

\bth
\label{cnm(1):th}
Fix $n\ge2$ and write
\beq
\label{cnm(1):eq}
c_{n,m}(1)=\sum_{k=0}^{n-1} a_k \binom{m}{k}
\eeq
for certain scalars $a_k$ depending on $n$.  We have the following.
\ben
\item[(a)] $a_0=0$.
\item[(b)] $a_1=1$.
\item[(c)] $a_2= \dil\binom{n}{\fl{n/2}}-2$.
\item[(d)] $a_{n-1} = 1$.
\item[(e)] $a_i\in\bbP$ for all $i\in[n-1]$.
\een
\eth
\bprf
(a)  Setting $m=0$ in~\eqref{cnm(1):eq} we obtain
$c_{n,0}(1)=a_0$.  But since $n\ge2$ there are no $w\in C(1)$ of length $n$ and maximum $0$.  So, $a_0=0$.

\medskip

(b)  Using part (a) and setting $m=1$ in~\eqref{cnm(1):eq} gives
$c_{n,1}(1)=a_1$.  There is a unique $w\in C(1)$ of length $n$ and maximum $1$, namely $w=1^n$.  This confirms the value of $a_1$.

\medskip

(c)  By~\cite[Corollary 6.3]{SW:cpm}, we have
$$
c_{n,2}(1) = \binom{n}{\fl{n/2}}.
$$
On the other hand, using~\eqref{cnm(1):eq} as well as parts (a) and (b) gives
$$
c_{n,2}(1) = \binom{2}{1}+a_2\binom{2}{2} = 2 + a_2.
$$
Equating the two expressions for $a_2$ finishes the proof.

\medskip

(d)  By Proposition~\ref{cnm(1):pro} we have that the coefficient of $m^{n-1}$ in $c_{n,m}(1)$ is $1/(n-1)!$.  But in~\eqref{cnm(1):eq}, the only term contributing to this power of $m$ is
$$
a_{n-1}\binom{m}{n-1} = \frac{a_{n-1}}{(n-1)!} m(m-1)\cdots (m-n+2).
$$
Thus $a_{n-1}=1$.

\medskip

(e)  By~\cite[equation (2)]{SW:cpm} we have
\beq
\label{cmn:exp}
c_{m,n}(1)=\sum_{\la\ptn n}  f^\la g_m^\la
\eeq
where $f^\la$ is the number of standard Young tableaux of shape $\la$ and $g_m^\la$ is a polynomial in $m$ which we will describe below.  Note that $f^\la$ is a positive integer for all $\la$.  So we will consider the binomial expansion of $g_m^\la$.

If $\la=(\la_1,\la_2,\ldots,\la_l)$ then let $\la'' = (\la_2,\ldots,\la_l)$ and suppose $\la''\ptn n''$.  Then $g_m^\la$ is the coefficient of $x^{m-2}$ in the power series expansion of
\beq
\label{fP:exp}
\frac{\sum_{\pi} x^{\des\pi}}{(1-x)^{n''+1}}
=\sum_\pi x^{\des\pi}\sum_{m\ge0}\binom{m+n''}{n''} x^m
\eeq
where $\pi$ runs over the linear extensions of the poset $\fP_{\la''}$ associated with $\la''$, and $\des\pi$ is the number of descents of $\pi$.  
Thus the coefficient of $x^{m-2}$ in the term of this expansion corresponding to $\pi$ is
\beq
\label{coeff}
\binom{m+n''-\des\pi-2}{n''}.
\eeq
We now break the proof into  two cases depending on whether $\la''$ is a column or not.

First consider the contribution of some $\la$ where $\la''$ is not a column, that is, $\la_2\ge2$.  It follows that $\fP_{\la''}$ has (at least) two integers $k,k'$ such that $k<k'$ in $\fP_{\la''}$ and $k<k'$ as integers.  So, the descending permutation $\pi=n''\ldots 21$ is not a linear extension of $\fP_{\la''}$ and any $\pi$ in~\eqref{fP:exp} satisfies $\des\pi\le n''-2$.  Letting $p:=n''-\des\pi-2$ we see that $0\le p \le n''-2$.  We can now write, using Vandermonde's convolution,
$$
\binom{m+p}{n''} = \sum_{i=2}^{n''} \binom{p}{n''-i} \binom{m}{i}.
$$
Since the coefficients $\binom{p}{n''-2}$ are nonnegative integers for $i\in[n-1]$, we are done with this case.

Now consider what happens when $\la''=(1^{n''})$ for some $0\le n''\le n-1$.  Then the poset $\fP_{\la''}$ is the chain $n''<n''-1<\ldots<1$.  So, the only possible linear extension is $\pi=n''\ldots 2 1$ with $\des\pi=n''-1$.  Using~\eqref{coeff} we see that the  contribution in~\eqref{fP:exp} to $g_m^\la$ for this $\la''$ is 
$\binom{m-1}{n''}$.  And in this case $\la=(n-n'',1^{n''})$ so that, by the hook formula, $f^\la=\binom{n-1}{n''}$.  Thus, by applying the binomial recursion twice, the total contribution to~\eqref{cmn:exp} of these $\la$
$$
\sum_{n''=0}^{n-1} \binom{n-1}{n''}\binom{m-1}{n''}
=\sum_{n''=1}^{n-1} \binom{n-2}{n''-1}\binom{m}{n''}.
$$
The coefficients $\binom{n-2}{n''-1}$ are positive integers for $n''\in[n-1]$.  So we can conclude that the $a_i$ are in $\bbP$, as desired.
\eprf

To better understand the numbers $c_{n,m}(u)$, we refine them according to the number of distinct letters used in a word. For $m,n \ge 1$, let 
\[\mathcal C_{n,m,k}(u) :=  \{w \in C(u) \mid |w| = n,\ \max w \le m, \text{ and there are $k$ distinct elements in $w$}\}\] 
and write 
\[
c_{n,m,k}(u) = |\mathcal C_{n,m,k}(u)|.\] 
Clearly, 
\beq
\label{cmn(u)Sum}
c_{n,m}(u) = \sum_{k = 1}^{\min(m,n)} c_{n,m,k}(u).
\eeq

Let $w$ be a word whose set of distinct elements is $S= \{s_1<s_1<\dots < s_k\}.$ 
The {\em standardization} of $w$ is the word $\std(w)$ 
obtained by replacing each element $a_j$ with its rank $j$ in the set $S$. Thus the standardization is $k$-packed as defined in the previous section.
Let  \[\mathcal B_{n,k} := \{ w \in C(1) \mid |w| = n \text{ and } w \text{ is $k$-packed}\}\] and \[b_{n,k} = |\mathcal B_{n,k}|.\]

\bth
Fix $n \ge 1$.  For all $m,k \ge 1$ we have
\[c_{n,m,k}(1) = b_{n,k}\binom{m-1}{k-1}.\]
\eth

   \bprf Let $m,k \ge 1$ be given. For a word $w \in \mathcal C_{n,m,k}(1) \subseteq C(1)$, we will show that $v=\std(w)\in\cB_{n,k}(1)$.
   Proposition~\ref{C(1)} (b) shows that the first row of $P(w)$ consists entirely of 1's, so in particular, $1 \in w.$ Hence $w$ contains a set $S\subseteq \{2,\dots, m\}$ of $k-1$ additional distinct elements.
   Because of the presence of the $1$'s, the alphabet of $v$  is exactly $[k]$. Since standardization is order-preserving, it also preserves the length of a longest weakly increasing subsequence. 
   Moreover, the presence of 1 in $w$ forces standardization to fix 1, so $v$ and $w$ have the same length of the longest weakly increasing subsequence ending at 1. By condition (c) of Proposition~\ref{C(1)}, $v \in C(1).$ 
   Therefore $v \in \mathcal B_{n,k}(1)$. 
   
   Conversely, given an $(k-1)$-element subset $S \subseteq\{2,\dots,m\}$ and a word $v \in \mathcal B_{n,k}$, the inverse relabeling replaces $2,\dots,k$ by the elements of $S$ in increasing order, producing a unique word $w$ of length $n$ with $k$ distinct letters from $[m].$ Since the condition of Proposition~\ref{C(1)} (c) is preserved under such relabelings, we have $w \in \mathcal C_{n,m,k}(1)$. 
   
   Thus, the words counted by $c_{n,m,k}(1)$ are in bijection with pairs $(S,v)$ where $S$ is an $(k-1)$-element subset of $\{2,\dots, m\}$ and $v\in\mathcal B_{n,k}$. Since there are $b_{n,k}$ choices for $v$, it follows that $c_{n,m,k}(1) = b_{n,k}\binom{m-1}{k-1}.$\eprf

Combining the previous theorem with~\eqref{cmn(u)Sum}, we immediately get the following result.

\bco
Let $n$ be fixed and suppose $m \ge n$. Then
\beq
\label{cnm(1):neweq}
c_{n,m}(1)=\sum_{k=1}^{n} b_{n,k} \binom{m-1}{k-1}.
\eeq
\eco

We now investigate the properties of the coefficients in this expansion.  Because we have a  combinatorial interpretation for the $b_{n,k}$, our proofs are simpler than for the coefficients $a_k$ in Theorem~\ref{cnm(1):th}.

\bth
Suppose $n$ is fixed and $m \ge n$,  write
\beq
\label{cnm(1):eq2}
c_{n,m}(1)=\sum_{k=1}^{n} b_k \binom{m-1}{k-1}
\eeq
where $b_k = b_{n,k}$.  We have the following.
\ben
\item[(a)] $b_k\in\bbP$ for all $k\in[n]$.
\item[(b)] $b_1 = 1$.
\item[(c)] $b_2 = \dil\binom{n}{\fl{n/2}}-1$
\item[(d)] $b_{n} = 1$.
\een
\eth

\bprf 
(a) Let $k \in [n]$. Consider the word $w = k(k-1)\dots 21^{n-k+1}.$ It is easy to see that $\text{lwi}(w) =\text{lwi}(w,1)$. By Proposition~\ref{C(1)}(c), this implies that $w \in C(1).$ Since $w$ has length $n$ and alphabet $[k]$, we have $w \in \mathcal B_{n,k}$. Hence, $b_k = b_{n,k} \in \mathbb{P}.$

\medskip

(b) There is a unique $1$-packed word of length $n$, namely $w=1^n$. Thus, $b_1 = b_{n,1} = 1.$

\medskip

(c) Notice that by definition, $c_{n,2}(1) = b_{n,1}+ b_{n,2}$. By~\cite[Corollary 6.3]{SW:cpm}, we have
\[c_{n,2}(1) = \binom{n}{\fl{n/2}}.\]
Using part (b) gives \[b_2 = b_{n,2} = \binom{n}{\fl{n/2}} - 1.\]

\medskip

(d) There is a unique $n$-packed word of length $n$, namely $w = n(n-1)\dots1$. Thus, $b_{n} =  b_{n,n} = 1.$\eprf 

We have a conjecture about the $b_{n,k}$ which we have verified by computer for $n \le 20$.

\begin{conj}
For fixed $n$, the sequence $\{b_{n,k}\}_{k\ge1}$ is log-concave.
\end{conj}

\section*{Acknowledgement} This material is based upon work supported by the National Science Foundation under Grant No. DMS-1929284 while the authors were in residence at the Institute for Computational and Experimental Research in Mathematics in Providence, RI, during the Categorification and Computation in Algebraic Combinatorics semester program.

\nocite{*}
\bibliographystyle{alpha}

\end{document}